\documentclass[sn-mathphys-num]{sn-jnl} 
\usepackage{graphicx}
\usepackage{multirow}
\usepackage{amsmath,amssymb,amsfonts}
\usepackage{amsthm}
\usepackage{mathrsfs}
\usepackage[title]{appendix}
\usepackage{xcolor}
\usepackage{textcomp}
\usepackage{manyfoot}
\usepackage{booktabs}
\usepackage{algorithm}
\usepackage{algorithmicx}
\usepackage{algpseudocode}
\usepackage{listings}

\theoremstyle{plain}
\newtheorem{theorem}{Theorem}
\newtheorem{proposition}[theorem]{Proposition}

\theoremstyle{remark}
\newtheorem{example}[theorem]{Example}

\theoremstyle{definition}
\newtheorem{definition}[theorem]{Definition}

\numberwithin{equation}{section}
\numberwithin{theorem}{section}

\usepackage{geometry}
\begin{document}

\title[Existence and Uniqueness of Mild Solutions to Rough Stochastic Partial Differential Equations driven by Fractional Brownian Motion with Hurst Index $H>1/4$]{\bf Existence and Uniqueness of Mild Solutions to Rough Stochastic Partial Differential Equations driven by Fractional Brownian Motion with Hurst Index $H>1/4$}

\author{\textbf{Maksim Vaskouski} \textbf{and} \textbf{Maksim Firsau}}

\affil{\it \orgdiv{Department of Fundamental Mathematics and Intelligent Systems}, \orgname{Belarusian State University}, \orgaddress{\street{4, Nezavisimosti Ave.}, \city{Minsk}, \postcode{220030}, \country{Belarus}} \par
{\centering \textit{Emails:} \texttt{vaskovskii@bsu.by; firsovm23@gmail.com}}}

\abstract{In this paper we prove a theorem on the unique solvability of the Cauchy problem for a stochastic partial differential equation driven by the fractional Brownian motion with Hurst index $H>\frac{1}{4}$.}


\maketitle

\section{Introduction}\label{sec1}

In the present paper we consider a  stochastic partial differential equation on a separable Hilbert space
\begin{equation}\label{eq1.1}
	dY_t=LY_tdt+F(Y_t)dB^H_t, \ t \in [0,T],
\end{equation}
with a negative self-adjoint differential operator $L$ generating an analytic semigroup $S_t=e^{Lt}$ and driven by a $d$-dimensional fractional Brownian motion $B^H_t$ with Hurst index $H > \frac{1}{4}.$ 

Stochastic differential equations driven by fractional Brownian motions have a lot of applications. Kleptsyna et al. \cite{bib1}
use models with multivariate
fractional Brownian motion to describe signal processes in filtering systems.
Cheridito \cite{bib2} considers the Samuelson model for stock price movements using
fractional Brownian motion. Zahle \cite{bib3} studies mixed type stochastic differential
equations in context of bond and stock models. Some other financial applications of
such equations reader can also find in Mishura's book \cite[Chapter 5]{bib4}.

Clearly, Eq. (\ref{eq1.1}) should be interpreted in the sense of the integral equation
\begin{equation}\label{eq1.2}
Y_t=S_tY_0+\int\limits_0^t S_{t-s}F(Y_s)dB^H_s, \ \ t\in [0,T].
\end{equation}
The main challenge here is how to understand the integral in the right-hand side of Eq.~(\ref{eq1.2}). Since $B^H_t$ is not a martingale for $H \ne \frac{1}{2}$, the classical Ito calculus is not applicable. Pathwise Young integrals can not be used due to lack of H\"older-regularity of $B^H_t$ with $H<\frac{1}{2}.$ Starting since $2000^{\textrm{th}}$ mathematicians apply successfully the rough path theory for such equations.  The original rough path theory was built by Lyons \cite{bib5} for ordinary differential equations driven by rough paths of arbitrary positive H\"older exponent, and this theory was applied by Coutin and Quan \cite{bib6} to derive existence and uniqueness results for Eq. (\ref{eq1.2}) in a finite-dimensional Euclidean space with $L=0$ driven by $B^H$, $H>\frac{1}{4}$. The threshold $\frac{1}{4}$ is very substantial: the second level of dyadic approximations has no a converging in probability sub-sequence in the case $H\le \frac{1}{4}.$ Although Nualart and Tindel \cite{bib7} have constructed a geometric rough path above $B^H$ for arbitrary $H>0,$ their construction is not compatible with the Lyons approach. Gubinelli \cite{bib8} has developed a functional version of the rough path theory based on notions of controlled rough path and rough path derivative, his theory covers $\alpha$-H\"older paths with $\alpha>\frac{1}{3}.$ Harang \cite{bib9} has extended the Gubinelli theory for alpha-H\"older rough paths with $\alpha>\frac{1}{4}.$ The theory of rough paths with an arbitrary H\"older exponent $\alpha$ for Eq. (\ref{eq1.2}) in an one-dimensional Euclidean space with $L = 0$ was developed in \cite{bib10}, \cite{bib11}, \cite{bib12}. Gerasimovics and Hairer \cite{bib13} have extended the Gubinelli theory for $\alpha$-H\"older rough paths, $\alpha>\frac{1}{3},$ taking values in a separable Hilbert space and controlled by an analytic semigroup. In the same work a theorem on the unique solvability of the Cauchy problem for equation (\ref{eq1.1}) is proved in the case $H > \frac{1}{3}$. In \cite{bib14} the theory of rough paths in infinite-dimensional spaces was extended to the case where rough paths have H\"older exponent $\alpha > \frac{1}{4}.$

In this paper, using the results of \cite{bib14}, we prove a theorem on the unique solvability of the Cauchy problem for equation (\ref{eq1.1}), assuming that the function $F$ is sufficiently smooth. We do not assume any global boundedness conditions on the right-hand side $F$, so our solution exists until an explosion time.

\section{Preliminaries}\label{sec2}

	\textbf{Definition of geometric rough paths.}
	Let $\mathcal{H}$ be a separable Hilbert space, $B$ be a Banach space.
	\par Let $\alpha \in [0, 1].$ Denote by $ C^{\alpha}([0, T], B)$ the set of $\alpha$-H\"older functions $f : [0, T] \rightarrow B$ such that
	$$ \|f\|_{\alpha, B} := \sup\limits_{s,t \in [0, T], s < t} \frac{\|f(t) - f(s)\|_{B}}{|t - s|^{\alpha}}
	< \infty,$$
	
	\par Let $ C_{2}^{\alpha}([0, T]^{2}, B)$ be the set of functions $g : [0, T] \times [0, T]  \rightarrow B$ such that
	$$ \|g\|_{\alpha, B} := \sup\limits_{s,t \in [0, T], s < t} \frac{\|g(s, t)\|_{B}}{|t - s|^{\alpha}}
	< \infty.$$
	
	\par For $f \in C^{\alpha}([0, T], B)$ we set $f_{s, t} := f(t) - f(s),$ $s,t \in [0, T],$ as an element of $C_{2}^{\alpha}([0, T]^{2}, B)$.
	
	\par Let $L^{\infty}([0, T], B)$ be the set of bounded functions $f: [0, T] \rightarrow B$ such that
	$$ 
	\|f\|_{\infty, B} := \sup\limits_{t \in [0, T]} \|f(t)\|_{B}<\infty.
	$$ 
	
	For arbitrary function $X: [0, T] \rightarrow \mathcal{H} $ an element $\textbf{X}^{1}_{s, t} := X_{s, t} \ \forall (s, t) \in [0, T]^{2},$ is called the first order process above $X$; elements $\textbf{X}^{2}: [0, T]^{2} \rightarrow \mathcal{H} \otimes \mathcal{H},$  $\textbf{X}^{3}: [0, T]^{2} \rightarrow \mathcal{H} \otimes \mathcal{H} \otimes \mathcal{H}$ are called the second and third order processes above $X$ if they satisfy to the following Chen's identity: 
	$$
	\textbf{X}^{2}_{s,t} - \textbf{X}^{2}_{s,u} - \textbf{X}^{2}_{u,t} = \textbf{X}^{1}_{s,u} \otimes \textbf{X}^{1}_{u,t}, 
	$$
	$$
	\textbf{X}^{3}_{s,t} - \textbf{X}^{3}_{s,u} - \textbf{X}^{3}_{u,t} = \textbf{X}^{1}_{s,u} \otimes \textbf{X}^{2}_{u,t} + \textbf{X}^{2}_{s,u} \otimes \textbf{X}^{1}_{u,t},
	$$
	for any $(s,u,t) \in [0, T]^{3}$.
	\par Let $\alpha \in (\frac{1}{4}, \frac{1}{3}]$. Denote by $\mathscr{C}^{\alpha}([0, T], \mathcal{H})$ the set of all
	$\alpha$-H\"older rough paths above $X \in C^{\alpha}([0, T], \mathcal{H})$ which are tuples $(1, \textbf{X}^{1}, \textbf{X}^{2}, \textbf{X}^{3}) $ such that $\textbf{X}^{1} \in C^{\alpha}_{2}([0, T]^{2}, \mathcal{H}),$ $\textbf{X}^{2} \in C^{2\alpha}_{2}([0, T]^{2}, \mathcal{H} \otimes \mathcal{H}),$ 
	$\textbf{X}^{3} \in C^{3\alpha}_{2}([0, T]^{2}, \mathcal{H} \otimes \mathcal{H} \otimes \mathcal{H})$ are processes of the first, second and third order above $X$.
	\par Let $ \mathscr{C}^{\alpha}_{g}([0, T], \mathcal{H})$ be the set of $\alpha$-H\"older geometric rough paths above $X \in C^{\alpha}([0, T], \mathcal{H})$ which are elements $(1, \textbf{X}^{1}, \textbf{X}^{2}, \textbf{X}^{3}) \in \mathscr{C}^{\alpha}([0, T], \mathcal{H})$ such that
	$$ Sym(\textbf{X}^{2}_{s, t}) = \frac{1}{2}\textbf{X}^{1}_{s, t} \otimes \textbf{X}^{1}_{s, t}, $$
	$$ Sym(\textbf{X}^{3}_{s, t}) = \frac{1}{6}\left(\textbf{X}^{1}_{s, t}\right)^{\otimes 3} $$
	for any $(s,t) \in [0, T]^{2}$.

	%\bigskip
		
	\textbf{Semigroups of linear operators.}
	Let $\mathcal{U}$ be a separable Hilbert space; $\mathcal{L}(\mathcal{U})$ be the set of bounded linear operators on $\mathcal{U}$; $\mathcal{L}_{2}(\mathcal{H}, \mathcal{U})$ be the set of Hilbert-Shmidt operators from $\mathcal{H}$ to $\mathcal{U}$.
	
	Recall that a family $\{S_{t}\}_{t\geq0}\subset\mathcal{L}(\mathcal{U})$ is called $C_0$-semigroup on $\mathcal{U}$ corresponding to a generator $A:D(A) \to \mathcal{U}$ if $D(A)$ is closed and dense in $\mathcal{U}$, $S_{t+s}=S_{t}S_{s} \ \forall s,t \geq 0,$ $S_{0} = I,$ $\|S_{t}u-S_{t_{0}}u\|_{\mathcal{U}}\underset{t\to t_{0}}{\longrightarrow} 0 \ \forall u \in \mathcal{U},$ and $Au=\displaystyle\lim_{t \to +0}\frac{S_{t}u-u}{t} \ \forall u \in D(A).$

	$C_0$-semigroup $\{S_{t}\}_{t\geq0}$ on $\mathcal{U}$ is called analytic if it can be extended on some sector
	$
	\Delta_{\delta} = \{\lambda \in \mathbb{C}: |arg\lambda| < \delta, Re\lambda > 0\}, \ 0 < \delta \leq \frac{\pi}{2},
	$
	such that $S_{\lambda}$ is continuous in $\overline{\Delta}_{\delta}$.

	Let $L$ be a negative self-adjoint linear operator on $\mathcal{U}$ with dense and closed domain $D(L) \subset \mathcal{U}$ such that
	$$
	\langle u, Lu\rangle_{\mathcal{U}} \leq C\langle u, u\rangle_{\mathcal{U}} \ \forall u \in D(L)
	$$
	for some constant $C<0.$ It is known that an operator $L$ satisfying the mentioned restrictions generates an analytic semigroup $\{S_{t}\}_{t\geq0}$ \cite[P. 421]{bib15}.

	For any $\beta \ge 0$ we define a separable Hilbert space $\mathcal{U}_{\beta}=D((-L)^{\beta})$ with the norm $\|\cdot\|_{\mathcal{U}_{\beta}} = \|(-L)^{\beta}\cdot\|_{\mathcal{U}}$. Also we define $\mathcal{U}_{-\beta}$ as the completion of $\mathcal{U}$ with respect to the norm $\|\cdot\|_{\mathcal{U}_{-\beta}} = \|(-L)^{-\beta}\cdot\|_{\mathcal{U}}$.

	Define the operator  $\tilde{S}_{t} = S_{t} - I$ and recall that for any $\rho, \beta, \gamma$ such that $\beta \geq \gamma,$  $\rho \in [0, 1]$ there hold estimates (see \cite[P. 41]{bib16}, \cite[P. 74]{bib17})
	$$
	\|S_{t}u\|_{\mathcal{U}_{\beta}} \lesssim t^{\gamma - \beta}\|u\|_{\mathcal{U}_{\gamma}}, \ \|\tilde{S}_{t}u\|_{\mathcal{U}_{\gamma - \rho}} \lesssim t^{\rho}\|u\|_{\mathcal{U}_{\gamma}},
	$$
	$\forall t \in (0, 1],$ $u \in \mathcal{U}_{\gamma}$, where $y \lesssim x$ denotes $y \leq Cx$ with some universal positive constant $C>0$ for arbitrary positive $x,$ $y$.
	
	%\bigskip
	
	\textbf{Integration of controlled rough paths.}
	Denote by $ C^{\alpha}_{S}([0, T], B)$ the set of functions $f : [0, T] \rightarrow B$ such that
	$ \|\delta_{S}f\|_{\alpha, B} < \infty,$
	where $\delta_{S}f_{s, t} = f_{t} - S_{t-s}f_{s}.$ Let $ C^{\alpha, 0}_{S}([0, T], B)$ be the set of functions $f : [0, T] \rightarrow B$ such that $ \|\delta_{S}f\|_{\alpha, 0, B} < \infty.$

	Let $\textbf{X} = (1, \textbf{X}^{1},\textbf{X}^{2},\textbf{X}^{3}) \in \mathscr{C}^{\alpha}_{g}([0, T], \mathcal{H})$, $\alpha \in (\frac{1}{4}, \frac{1}{3}]$. We shall say that a function $Y \in C^{\alpha}_{S}([0, T], B)$ is controlling by the geometric rough path $\textbf{X}$ with respect to the semigroup $\{S_{t}\}_{t\geq0}$ if there exist elements $Y' \in C^{\alpha}_{S}([0, T], \mathcal{L}_{2}(\mathcal{H}, B)),$  $Y'' \in C^{\alpha}_{S}([0, T], \mathcal{L}_{2}(\mathcal{H}, \mathcal{L}_{2}(\mathcal{H}, B)))$ such that $ \|R^{S,Y,2}\|_{2\alpha, \mathcal{L}_{2}(\mathcal{H}, B)} < \infty,$ $ \|R^{S,Y,3}\|_{3\alpha, B} < \infty$, where 
	$$ 
	R^{S,Y,2}_{s, t} := Y_{t}' - S_{t-s}(Y_{s}' + Y_{s}''\textbf{X}^{1}_{s, t}),
	$$
	$$
	R^{S,Y,3}_{s, t} := Y_{t} - S_{t-s}(Y_{s} + Y_{s}'\textbf{X}^{1}_{s, t} + Y_{s}''\textbf{X}^{2}_{s, t}),
	$$
	and the element $(Y,Y',Y'')$ is called $\alpha$-H\"older controlled rough path with respect to $\textbf{X}$ and $\{S_{t}\}_{t\geq0}$. In further we will write $\|\cdot\|_{B}$ instead of $\|\cdot\|_{\mathcal{L}_{2}(\mathcal{H}, B)},$ $\|\cdot\|_{\mathcal{L}_{2}(\mathcal{H}, \mathcal{L}_{2}(\mathcal{H}, B))}$ and we will write $\|\cdot\|_{\alpha, B}$ instead of $\|\cdot\|_{\alpha, \mathcal{L}_{2}(\mathcal{H}, B)},$ $\|\cdot\|_{\alpha, \mathcal{L}_{2}(\mathcal{H}, \mathcal{L}_{2}(\mathcal{H}, B))}$.
	
	We define a Banach space of $\alpha$-H\"older controlled rough path with respect to $\textbf{X}$ and $\{S_{t}\}_{t\geq0}$:
	$$ 
	\mathcal{D}^{\alpha}_{S,\textbf{X}}([0, T], B) := \biggl\{(Y, Y', Y''): Y \in C^{\alpha}_{S}([0, T], B),
	$$
	$$
	Y' \in C^{\alpha}_{S}([0, T], \mathcal{L}_{2}(\mathcal{H}, B)), \sum_{i=1}^{3}\|R^{S,Y,i}\|_{i\alpha, B} < \infty\biggr\}
	$$ 
with the norm	
	$$
	\|(Y, Y', Y'')\|_{\mathcal{D}^{\alpha}_{S,\textbf{X}}} := \sum_{i=0}^{2}\|Y^{(i)}_{0}\|_{B} + \|(Y, Y', Y'')\|_{S, \textbf{X}, \alpha, B},
	$$
	where $\|(Y, Y', Y'')\|_{S, \textbf{X}, \alpha, B} := \sum\limits_{i=1}^{3}\|R^{S,Y,i}\|_{i\alpha, B}, R^{S,Y,1}_{s, t} :=  Y_{t}'' - S_{t-s}Y_{s}''.$

	Let $s,t \in [0, T], s < t, \alpha \in (\frac{1}{4}, \frac{1}{3}], \textbf{X} = (1, \textbf{X}^{1},\textbf{X}^{2},\textbf{X}^{3}) \in \mathscr{C}^{\alpha}([0, T], \mathcal{H})$ and $(Y, Y', Y'') \in \mathcal{D}^{\alpha}_{S,\textbf{X}}([0, T], \mathcal{L}_{2}(\mathcal{H}, \mathcal{U}_{\beta}))$. We define rough path integral of the convolution $S_{t-u}Y_u$ over the rough path $\textbf{X}$ as the following limit
	\begin{equation}\label{eq2.1}
		\int\limits_s^t S_{t - u}Y_{u}d\textbf{X}_{u} := \lim_{|\mathcal{P}| \to 0}\sum_{[u,v] \in \mathcal{P}} S_{t - u}(Y_{u}\textbf{X}^{1}_{u, v} + Y_{u}'\textbf{X}^{2}_{u, v} + Y_{u}''\textbf{X}^{3}_{u, v}),
	\end{equation}
	where $\mathcal{P}$ is a finite partition of $[s, t],$  $|\mathcal{P}| = \max\limits_{[u,v] \in \mathcal{P}}|v - u|,$ assuming that the limit is finite and does not depend on partitions of segment $[s,t].$
	
\textbf{Fractional Brownian Motion.}
Let $H \in (0,1).$ A scalar fractional Brownian motion $B^{H} = \{B^{H}_{t}: t \in \mathbb{R}^+\}$ is a Gaussian process defined on a complete probability space $(\Omega, \mathcal{F}, \mathbb{P})$ such that
	\begin{enumerate} 
		\item $B^{H}_{0} = 0$ a.s.;
		\item $\mathbb{E}(B^{H}_{t}) = 0, t \in \mathbb{R}^+;$
		\item $\mathbb{E}(B^{H}_{t}B^{H}_{s}) = \frac{1}{2}(|t|^{2H} + |s|^{2H} - |t - s|^{2H}), s, t \in \mathbb{R}^+.$
	\end{enumerate}
	
	Let $H \in (0,1),$ $B_i^H,$ $i=1,\ldots,d$ be scalar independent fractional Brownian motions, then the process $B^H=(B_1^H,\ldots, B_d^H)$ is called $d$-dimensional fractional Brownian motion.

	\begin{proposition} (Theorem 1.1, \cite{bib7}) \label{th2.4}
	Let $H \in (\frac{1}{4}, \frac{1}{3}]$ and $B^H$ be a $d$-dimensional fractional Brownian motion with Hurst index $H.$ Then for any $\alpha \in (\frac{1}{4}, H)$ there exists a stochastic process $\mathbb{B},$ adapted to the filtration generated by $B^H,$ such that $\mathbb{B} = (1, \mathbb{B}^{1}, \mathbb{B}^{2}, \mathbb{B}^{3}) \in \mathscr{C}^{\alpha}_{g}([0, T], \mathbb{R}^{d})$ is  $\alpha$-H\"older geometric rough path above $B^H$ almost surely.
	\end{proposition}
	
	\textbf{One-point extension of a separable Hilbert space}. Let $X$ be a nonempty metric space with metric  $d.$ Any element $x \in X$ can be represented as the equivalence class of all sequences $(x_{n}) \subset X$ such that
	$$
		d(x_{n}, x)  \underset{n \rightarrow \infty}{\longrightarrow} 0.
	$$
	Denote by \mbox{\mathversion{bold}$\infty$} the equivalence class of all sequences $(x_{n}) \subset X$ such that
	$$
	d(0, x_{n})  \underset{n \rightarrow \infty}{\longrightarrow} \infty. 
	$$
	Then the space $\widehat{X}$ is called one-point extension of $X$, if $\widehat{X} = X \cup \{\mbox{\mathversion{bold}$\infty$}\}$ and $\widehat{X}$ is metrizable.
	
	\begin{proposition}\label{th2.5}
	Any separable Hilbert space $\mathcal{U}$ admits one-point extension $\widehat{\mathcal{U}} =\mathcal{U} \cup \{\mbox{\mathversion{bold}$\infty$}\}$ with metric 
	\begin{equation}\label{eq2.2}
	\widehat{d}(u, v) = \left\{\begin{matrix}
	\frac{\|u - v\|}{\sqrt{(1 + \|u\|^{2})(1 + \|v\|^{2})}}, & u \in \mathcal{U}, v \in \mathcal{U},
	\\ \frac{1}{\sqrt{1 + \|u\|^{2}}}, & u \in \mathcal{U}, v = \mbox{\mathversion{bold}$\infty$},
	\\ 0, & u = \mbox{\mathversion{bold}$\infty$}, v = \mbox{\mathversion{bold}$\infty$},
	\end{matrix}\right.
	\end{equation}
	where $\|\cdot\|$ is the norm of $\mathcal{U}$ induced by the inner product. Moreover, $\widehat{\mathcal{U}}$ is complete with respect to the metric $\widehat{d}.$
	\end{proposition}
	See Appendix A.1 for a proof.

	\section{Existence and uniqueness of mild solutions for rough stochastic PDEs driven by fBm}\label{sec3}
	
	Let $B^{H}$ be a $d$-dimensional fractional Brownian motion with Hurst index $H \in (\frac{1}{4}, \frac{1}{3}]$ defined on a complete probability space $(\Omega, \mathcal{F}, \mathbb{P}),$ $\mathbb{B} = (1, \mathbb{B}^{1}, \mathbb{B}^{2}, \mathbb{B}^{3}) \in \mathscr{C}^{\alpha}_{g}([0, T], \mathbb{R}^{d})$ be $\alpha$-H\"older geometric rough path above $B^H,$  $\alpha \in (\frac{1}{4}, H),$ $\xi$ be a random variable on $(\Omega, \mathcal{F}, \mathbb{P})$ taking values in a separable Hilbert space $\mathcal{U},$ $F:\mathcal{U} \rightarrow \mathcal{L}_{2}(\mathbb{R}^{d},\mathcal{U})$ and $L$ be a negative self-adjoint linear operator generating an analytic semigroup $\{S_t\}_{t \ge 0}$ on $\mathcal{U}$.

	We consider a rough path evolution equation
	\begin{equation}\label{eq5.1}
	dY_{t} = LY_{t}dt + F(Y_{t})d\mathbb{B}_{t},\quad
	Y_{0} = \xi.
	\end{equation}
	
	Let $\widehat{\mathcal{U}} =\mathcal{U} \cup \{\mbox{\mathversion{bold}$\infty$}\}$ -- be the one-point extension of $\mathcal{U}$ with metric (\ref{eq2.2}). For arbitrary function $f:[0,T]\to \widehat{\mathcal{U}}$ value $\tau(f) =\inf\{t: f(t) = \mbox{\mathversion{bold}$\infty$}\}$ is called an explosion time of $f.$ If $\{t: f(t) = \mbox{\mathversion{bold}$\infty$}\} = \varnothing,$ we set $\tau(f) = +\infty.$ 
	
	Denote by $\widehat{C}([0, T], \widehat{\mathcal{U}})$ the set of functions $f:[0,T]\to \widehat{\mathcal{U}}$ such that
 	$f$ is continuous for $t \in [0, T] \cap [0,\tau(f)),$ and $f(t) = \mbox{\mathversion{bold}$\infty$}$ for $t \in [0, T] \cap[\tau(f),\infty).$ 

	Define the Banach space of controlled rough paths
	$$
	\tilde{\mathcal{D}}^{\alpha}_{S,\textbf{X}}([0, T],\mathcal{U}_{\beta}) := \mathcal{D}^{\alpha}_{S,\textbf{X}}([0, T],\mathcal{U}_{\beta - 3\alpha}) \cap ( C^{\alpha}_{S}([0, T],\mathcal{U}_{\beta}) \times  L^{\infty}([0, T], \mathcal{L}_{2}(\mathcal{H},\mathcal{U}_{\beta}))^{\times 2})
	$$
	with the norm
	$$
	\|(Y, Y', Y'')\|_{\tilde{\mathcal{D}}^{\alpha}_{S,\textbf{X}}([0, T],\mathcal{U}_{\beta})} := \|\delta_{S}Y\|_{\alpha,\mathcal{U}_{\beta}} + \|Y'\|_{\infty,\mathcal{U}_{\beta}} + \|Y''\|_{\infty,\mathcal{U}_{\beta}} 
	$$
	$$
	+ \|(Y, Y', Y'')\|_{\mathcal{D}^{\alpha}_{S,\textbf{X}}([0, T],\mathcal{U}_{\beta - 3\alpha})}.
	$$

	Further we will write $L^{\infty}([0, T], B)^{\times 2}$ instead of $L^{\infty}([0, T], B) \times L^{\infty}([0, T], \mathcal{L}_{2}(\mathcal{H}, B)).$ The same applies to $L^{\infty}([0, T], B)^{\times n}$ for $n > 2.$

	\begin{definition}\label{th5.1}
	 Let $H \in (\frac{1}{4},\frac{1}{3}],$ $\alpha \in (\frac{1}{4},H),$ $T > 0$. A stochastic process $Y_{t}$ defined on $(\Omega, \mathcal{F}, \mathbb{P})$ is called mild solution of Caushy problem (\ref{eq5.1}) on segment $[0, T]$ if there hold the following conditions:
	
	1) $Y(\omega) \in \widehat{C}([0, T], \widehat{\mathcal{U}})$ a.s.;
	
	2) for almost all $\omega \in \Omega$  and any $\tilde{T} \in (0, \tau(\omega)) \cap [0, T]$ the mapping $t \rightarrow Y_{t}(\omega)$ is a mild solution of the deterministic equation
	$$
	dY_{t} = LY_{t}dt + F(Y_{t})d\mathbb{B}_{t}(\omega),\quad t \in [0, \tilde{T}],
	$$
	with initial condition $Y_{0}(\omega) = \xi(\omega),$ 	where $\tau(\omega) = \tau(Y(\omega))$ is the explosion time of $Y(\omega),$ that is:
		
	\quad 2.1) $(Y, Y', Y'')(\omega) \in \tilde{\mathcal{D}}^{\alpha}_{S,\textbf{X}}([0, \tilde{T}],\mathcal{U})$, where $Y' = F(Y),$ $Y'' = DF(Y)F(Y)$;
	
	\quad 2.2) $(F(Y), F(Y)',F(Y)'')(\omega) \in \mathcal{D}^{\alpha}_{S,\textbf{X}}([0, \tilde{T}], \mathcal{L}_{2}(\mathcal{H},\mathcal{U}_{-3\alpha}))$ \\
	$\cap$ $L^{\infty}([0, \tilde{T}], \mathcal{L}_{2}(\mathcal{H},\mathcal{U}))^{\times 3};$
	
	\quad 2.3) for any $t \in [0, \tilde{T}]$ there holds
	$$
	Y_{t}(\omega) = S_{t}\xi(\omega) + \int_{0}^{t}S_{t-s}F(Y_{s}(\omega))d\mathbb{B}_{s}(\omega),
	$$
	where the integral in the righ-hand side is the rough path integral $(\ref{eq2.1})$.
	
	3) for almost all $\omega \in \Omega$ there holds $\limsup\limits_{t \uparrow \tau(\omega)}\Vert Y_t(\omega) \Vert_\mathcal{U}=\infty$ assuming $\tau(\omega) \le T.$
	
	Mild solution $Y_{t}$ of Caushy problem (\ref{eq5.1}) on segment $[0, T]$ is called unique if there holds $\mathbb{P}(Y_{t} = \widetilde{Y}_{t} \  \forall t \in [0, T]) = 1$ for any other mild solution $\widetilde{Y}_{t}$ of Caushy problem.
	\end{definition}
	
	\begin{theorem}\label{th5.2}
	Let $H \in (\frac{1}{4},\frac{1}{3}],$ $\alpha \in (\frac{1}{4},H),$ $T > 0,$ $F \in C^{4}(\mathcal{U}_{\gamma},\mathcal{L}_{2}(\mathbb{R}^{d},\mathcal{U}_{\gamma}))$ for any $\gamma \in (-\alpha, 0]$. Then there exists a unique mild solution $Y_{t}$ of Caushy problem (\ref{eq5.1}) on segment $[0, T].$
	\end{theorem}

\begin{proof}
	According to Theorem \cite{bib14} for almost all $\omega \in \Omega$ there exists $T_{0}(\omega) \in (0,T]$ such that deterministic equation (\ref{eq5.1}) with initial condition $Y_{0}(\omega) = \xi(\omega)$ has a unique mild solution $t \rightarrow Y_{t}(\omega)$ on segment $[0, T_{0}(\omega)]$. We have
	$$
	Y_{t}(\omega) = S_{t}Y_{0}(\omega) + \int_{0}^{t}S_{t-s}F(Y_{s}(\omega))d\mathbb{B}_{s}(\omega), \ t \in [0, T_{0}(\omega)].
	$$
	If $T_0(\omega) < T$, we extend the constructed solution $Y_t(\omega)$ to the right of the point $T_0(\omega)$ as follows.
	
	For almost all $\omega$ the element $\tilde{\mathbb{B}}_{\cdot}(\omega) = \mathbb{B}_{\cdot + T_0(\omega)}(\omega)$ belongs to $\mathscr{C}^{\alpha}_{g}([0, T - T_0(\omega)], \mathbb{R}^{d})$. By Theorem \cite{bib14} we will construct a deterministic solution $Z_t(\omega)$, $t \in [0, T_1(\omega)]$, to the equation
	$$
		dZ_{t}(\omega) = LZ_{t}(\omega)dt + F(Z_{t}(\omega))d\tilde{\mathbb{B}}_{t}(\omega)
	$$
	with the initial condition $Z_{0}(\omega) = Y_{T_0(\omega)}(\omega)$, where $T_0(\omega) + T_1(\omega) \leqslant T$, then
	$$
	Z_{t}(\omega) = S_{t}Z_{0}(\omega) + \int_{0}^{t}S_{t-s}F(Z_{s}(\omega))d\tilde{\mathbb{B}}_{s}(\omega), \ t \in [0, T_{1}(\omega)].
	$$
	Define the mapping $\tilde{Y}_t(\omega)$, $t \in [0, T_0(\omega) + T_1(\omega)]$, as follows:
	$$
	\tilde{Y}_{t}(\omega) = Y_{t}(\omega), \ t\in [0,T_{0}(\omega )], \quad
	\tilde{Y}_{t}(\omega) = Z_{t-T_{0}(\omega)}(\omega), \ t\in (T_{0}(\omega ),T_{0}(\omega )+T_{1}(\omega)].
	$$
	Then, for every $t \in [0, T_0(\omega)]$, the relation holds
	$$
	\tilde{Y}_{t}(\omega) = S_{t}\tilde{Y}_{0}(\omega) + \int_{0}^{t}S_{t-s}F(\tilde{Y}_{s}(\omega))d\mathbb{B}_{s}(\omega).
	$$
	For any $t \in (T_0(\omega), T_0(\omega)+T_1(\omega)]$, we have
	$$
	\tilde{Y}_{t}(\omega) = Z_{t-T_{0}(\omega)}(\omega) = S_{t-T_{0}(\omega)}Z_{0}(\omega) + \int_{0}^{t-T_{0}(\omega)}S_{t-T_{0}(\omega)-s}F(Z_{s}(\omega))d\tilde{\mathbb{B}}_{s}(\omega)
	$$
	$$
	= S_{t-T_{0}(\omega)}Y_{T_0(\omega)}(\omega) + \int_{0}^{t-T_{0}(\omega)}S_{t-T_{0}(\omega)-s}F(\tilde{Y}_{s + T_{0}(\omega)}(\omega))d\mathbb{B}_{s + T_{0}(\omega)}(\omega)
	$$
	$$
	= S_{t-T_{0}(\omega)}\tilde{Y}_{T_0(\omega)}(\omega) + \int_{T_{0}(\omega)}^{t}S_{t-u}F(\tilde{Y}_{u}(\omega))d\mathbb{B}_{u}(\omega)
	$$
	$$
	= S_{t-T_{0}(\omega)}S_{T_0(\omega)}\tilde{Y}_{0}(\omega) + \int_{0}^{T_0(\omega)}S_{t-T_{0}(\omega)}S_{T_0(\omega)-s}F(\tilde{Y}_{s}(\omega))d\mathbb{B}_{s}(\omega)
	$$
	$$
	+ \int_{T_{0}(\omega)}^{t}S_{t-u}F(\tilde{Y}_{u}(\omega))d\mathbb{B}_{u}(\omega) = S_{t}\tilde{Y}_{0}(\omega) + \int_{0}^{t}S_{t-s}F(\tilde{Y}_{s}(\omega))d\mathbb{B}_{s}(\omega).
	$$
	
	Thus, $\tilde{Y}_t(\omega)$ is a deterministic mild solution to equation (\ref{eq5.1}) on $[0, T_0(\omega) + T_1(\omega)]$ with the initial condition $\tilde{Y}_0(\omega) = \xi(\omega)$. Instead of $\tilde{Y}_{t}(\omega)$, we will again write $Y_{t}(\omega)$.
	
	Applying the same procedure for any $n \in \mathbb{N},$ we will obtain a deterministic mild solution of equation (\ref{eq5.1}) on segment $[0, \tau_{0}^{n}(\omega)],$ where $\displaystyle \tau_{0}^{n}(\omega) := \sum_{k = 0}^{n}T_{k}(\omega).$ If for some natural $n_0 < n$ we obtain that $\displaystyle \sum_{k=0}^{n_0} T_k(\omega) \ge T$, then we set $T_k(\omega) = 0$ for every $k \in [n_0 + 1, n]$.
	
	Hence, for a fixed $\omega$ there are possible alternatives: either there exists $n \in \mathbb{N}$ such that $\tau_{0}^{n}(\omega) \geq T,$ or  $\tau_{0}^{n}(\omega) < T$ for any $n \in \mathbb{N}$ implying that the series $\displaystyle \sum_{k = 0}^{\infty}T_{k}(\omega)$ converges and $T_{n}(\omega) \underset{n \rightarrow \infty}{\longrightarrow} 0.$ 
	
	Let's consider the second alternative. We fix some $\omega \in \Omega.$ Let $\displaystyle \tau(\omega) := \sum_{k = 0}^{\infty}T_{k}(\omega).$ It's clear that $\tau_{0}^{n}(\omega) \underset{n \rightarrow \infty}{\longrightarrow} \tau(\omega)$ and $\tau(\omega) \leq T.$ By the way we obtain a unique deterministic mild solution $Y_{t}(\omega)$ of Caushy problem (\ref{eq5.1}) on $[0, \tau(\omega)).$ Since $\{S_{t}\}_{t\geq0}$ is a $C_0$-semigroup, we can show that $Y_{t}(\omega)$ is continuous on $t \in [0,\tau(\omega)):$
	$$
	\|Y_{t}(\omega) - Y_{t_{0}}(\omega)\|_{\mathcal{U}} \leq 	\|Y_{t}(\omega) - S_{t-t_{0}}Y_{t_{0}}(\omega)\|_{\mathcal{U}} + \|S_{t-t_{0}}Y_{t_{0}}(\omega) - S_{0}Y_{t_{0}}(\omega)\|_{\mathcal{U}}
	$$
	$$
	\leq \|\delta_{S}Y(\omega)\|_{\alpha,\mathcal{U}}|t-t_{0}|^{\alpha} + \|S_{t-t_{0}}Y_{t_{0}}(\omega) - S_{0}Y_{t_{0}}(\omega)\|_{\mathcal{U}} \underset{t \rightarrow t_{0}}{\longrightarrow} 0 \ \forall t_{0} \in [0,\tau(\omega)).
	$$
	Let's prove that $\limsup\limits_{t \uparrow \tau(\omega)}\|Y_{t}(\omega)\|_{\mathcal{U}} = \infty.$
	
	It follows from the proof of Theorem \cite{bib14} that for any $\tilde{\varepsilon} \in (0, 1)$ and any $n \in \mathbb{N} \cup \{0\}$ value $T_{n}(\omega)$ can be chosen as follows:
	\begin{equation}\label{eq5.2}
	T_{n}(\omega) = \frac{1 - \tilde{\varepsilon}}{C_{\alpha, \varepsilon, \mathbb{B}(\omega)}(1 + \|F\|_{C^{4}_{b}(E_{3\varepsilon-3\alpha}^{n}(\omega),\mathcal{U}_{3\varepsilon-3\alpha})})^{p}(1 + \|\xi^{n}(\omega)\|_{\mathcal{U}})^{q}},
	\end{equation}
	where $p, q > 0$ are some constants which depend only on $\alpha,$ $\varepsilon,$ $\xi^{n}(\omega) = Y_{\tau_{0}^{n-1}(\omega)}(\omega),$ $n \geq 1,$ $\xi^{0}(\omega) = \xi(\omega),$ and $(E_{3\varepsilon-3\alpha}^{n}(\omega))$ is a sequence of balls in $\mathcal{U}_{3\varepsilon-3\alpha}$ such that 
	$$
	E_{3\varepsilon-3\alpha}^{n}(\omega) = \{u \in\mathcal{U}_{3\varepsilon - 3\alpha}: \|u\|_{\mathcal{U}_{3\varepsilon - 3\alpha}} \leq 1 + \|\bar{Y}^{n}(\omega)\|_{\infty,\mathcal{U}_{3\varepsilon - 3\alpha}}\}.
	$$
	
	From the proof of Theorem \cite{bib14}, it follows that the constant $C_{\alpha,\varepsilon,\mathbb{B}(\omega)}$ can be chosen such that the map $\omega \mapsto C_{\alpha,\varepsilon,\mathbb{B}(\omega)}$ is a random variable. Consequently, $T_0(\omega)$ is a random variable. Furthermore, it follows from the proof of Theorem \cite{bib14} that for any random variable $\eta(\omega)$,  such that $0 < \eta(\omega) \leq T_0(\omega)$ almost surely, the map $\omega \mapsto Y_{\eta(\omega)}(\omega)$ is a random variable. By induction, we obtain that for every natural number $n$ the map $\omega \mapsto T_n(\omega)$ is a random variable. It follows that $\tau(\omega)$ is also a random variable.
	
	The constant $C_{\alpha, \varepsilon, \mathbb{B}(\omega)}$ does not depend on $n$ since the norm $\|\mathbb{B}^{i}(\omega)\|_{i\alpha, (\mathbb{R}^{d})^{\otimes i}}$ is the norm of the process $\mathbb{B}^{i}(\omega)$ defined on segment $[0, T], \forall i = \overline{1, 3}.$
	It's clear that
	$$
	\|\bar{Y}^{n}(\omega)\|_{\infty,\mathcal{U}_{3\varepsilon - 3\alpha}} \lesssim \|\delta_{S}\bar{Y}^{n}(\omega)\|_{\alpha,\mathcal{U}_{3\varepsilon - 3\alpha}} + \|\xi^{n}(\omega)\|_{\mathcal{U}}
	$$
	$$
	\lesssim \|F(\xi^{n}(\omega))\|_{\mathcal{U}}(\|\mathbb{B}^{1}(\omega)\|_{\alpha, \mathbb{R}^{d}} + \|DF(\xi^{n}(\omega))\|_{\mathcal{U}}\|\mathbb{B}^{2}(\omega)\|_{2\alpha, \mathbb{R}^{d} \otimes \mathbb{R}^{d}}) + \|\xi^{n}(\omega)\|_{\mathcal{U}}.
	$$
	Since the function $F$ is locally bounded on $\mathcal{U},$ there exists a strictly increasing function $f: \mathbb{R}^{+} \rightarrow \mathbb{R}^{+}$ such that $\|F(u)\|_{\mathcal{U}} \leq f(\|u\|_{\mathcal{U}}) \ \forall u \in\mathcal{U}.$ For example, we can take the following function $f$:
	$\displaystyle f(x) = \sup_{\tilde{u} \in E(x)}\|F(\tilde{u})\|_{\mathcal{U}},$ where $E(x) = \{\tilde{u} \in \mathcal{U}: \|\tilde{u}\|_{\mathcal{U}} \leq x\},  \ \forall x \in \mathbb{R}^{+}.$ The same will be true for $DF:$  $\|DF(u)\|_{\mathcal{U}} \leq g(\|u\|_{\mathcal{U}}) \ \forall u \in \mathcal{U},$ with some strictly increasing function $g$. Let $h(x) = 1 + x + f(x)(1 + g(x)).$ It's clear that there holds the estimate
	$$
	 \|u\|_{\mathcal{U}_{3\varepsilon - 3\alpha}} \lesssim h(\|\xi^{n}(\omega)\|_{\mathcal{U}})
	$$
	for any $u \in E_{3\varepsilon-3\alpha}^{n}(\omega).$
	Without loss of generality we may assume that for a fixed $n \in \mathbb{N}$ the function $F$ in (\ref{eq5.2}) is defined on the ball $\tilde{E}_{3\varepsilon-3\alpha}(\|\xi^{n}(\omega)\|_{\mathcal{U}}),$ where
	$$
	\tilde{E}_{3\varepsilon-3\alpha}(x) = \big\{u \in\mathcal{U}_{3\varepsilon - 3\alpha}: \|u\|_{\mathcal{U}_{3\varepsilon - 3\alpha}} \lesssim h(x)\big\}.
	$$
	Since the function $h(x)$ is strictly increasing, $\|F\|_{C^{4}_{b}(\tilde{E}_{3\varepsilon-3\alpha}(x),\mathcal{U}_{3\varepsilon-3\alpha})}$ is also an increasing function. Hence,
	$C_{\alpha, \varepsilon, \mathbb{B}(\omega)}(1 + \|F\|_{C^{4}_{b}(\tilde{E}_{3\varepsilon-3\alpha}(x),\mathcal{U}_{3\varepsilon-3\alpha})})^{p}(1 + x)^{q}$ is a strictly increasing function. Since $T_{n}(\omega) \underset{n \rightarrow \infty}{\longrightarrow} 0,$ we have $\|\xi^{n}(\omega)\|_{\mathcal{U}} \underset{n \rightarrow \infty}{\longrightarrow} \infty$ that implies $\|Y_{\tau_{0}^{n}(\omega)}(\omega)\|_{\mathcal{U}} \underset{n \rightarrow \infty}{\longrightarrow} \infty.$ So, $\displaystyle \limsup_{t \rightarrow \tau(\omega) - 0}\|Y_{t}(\omega)\|_{\mathcal{U}} = +\infty.$
	
	Let $Y_t(\omega) = \mbox{\mathversion{bold}$\infty$}$ for $t\in [\tau(\omega),T]$ if $\tau(\omega)\le T.$
	
	Clearly, $Y$ is a stochastic process whose paths belong to  $\widehat{C}([0, T], \widehat{\mathcal{U}})$ almost surely. The uniqueness of the constructed solution follows from Theorem \cite{bib14}.
	
	Theorem is proved.
\end{proof}
	
Applying the obtained Theorem, one can see that there exists a unique mild solution for the following two problems.
	
\begin{example}\label{th5.3}
	The Caushy problem for the stochastic heat equation.
	\begin{center}
		$
		\begin{cases}
			\frac{\partial u}{\partial t}(t, \varkappa) = \Delta u(t, \varkappa) + f(u(t, \varkappa)) \dot{B}^{H}_{t}, & t \in (0, T], \quad \varkappa \in \mathbb{T}^{3},\\
			u(0, \varkappa) = \varphi(\varkappa), & \varkappa \in \mathbb{T}^{3}.
		\end{cases}
		$
		
	\end{center}

	Here $Y_{t} = u(t, \cdot)$ is a stochastic processes taking values in the separable Hilbert space $\mathcal{U} = H^{k}(\mathbb{T}^{3}, \mathbb{R}),$ $Y_{0} = \varphi(\cdot),$ $H \in (\frac{1}{4}, \frac{1}{3}],$ $\alpha \in (\frac{1}{4}, H),$ $\mathbb{B} = (1, \mathbb{B}^{1}, \mathbb{B}^{2}, \mathbb{B}^{3}) \in \mathscr{C}^{\alpha}_{g}([0, T], \mathbb{R})$ -- a geometric rough path above $B^{H},$ $L = \Delta$ -- negative self-adjoint operator generating an analytic semigroup $\{e^{\Delta t}\}_{t\geq0}$ on $\mathcal{U}$ (see \cite[Ex. 12.17]{bib18}, \cite[Sec. 6.4]{bib16}), $D(\Delta) = H^{k + 2}(\mathbb{T}^{3}, \mathbb{R}),$ $\mathcal{U}_{\beta} = H^{k + 2\beta}(\mathbb{T}^{3}, \mathbb{R}),$ $k > \frac{3}{2} + 2\alpha,$ $\beta \in (-\alpha, 0],$ $f \in \mathbb{R}[x],$ $(F(y)z)(\varkappa) = f(y(\varkappa))z,$ $y \in\mathcal{U},$ $z \in \mathbb{R},$ $\varkappa \in \mathbb{T}^{3},$ $\mathbb{T}^{3}$ -- a three-dimensional torus. 
\end{example}
	
\begin{example}\label{th5.4}
	The initial-boundary problem for the stochastic heat equation.
	\begin{center}
		$
		\begin{cases}
			\frac{\partial u}{\partial t}(t, x) = 	\frac{\partial^2 u}{\partial x^2}(t, x) + u(t, x) \dot{B}^{H}_{t}, & t \in (0, T], \quad x \in (0, 1),\\
			u(0, x) = \varphi(x), & x \in [0, 1],
			\\ u(t, 0) = u(t, 1) = 0, & t \in [0, T].
		\end{cases}
		$
		
	\end{center}

	Here $Y_{t} = u(t, \cdot)$ is a stochastic processes taking values in the separable Hilbert space $\mathcal{U} = L^{2}([0, 1], \mathbb{R}),$ $Y_{0} = \varphi(\cdot),$ $H \in (\frac{1}{4}, \frac{1}{3}],$ $\alpha \in (\frac{1}{4}, H),$ $\mathbb{B} = (1, \mathbb{B}^{1}, \mathbb{B}^{2}, \mathbb{B}^{3}) \in \mathscr{C}^{\alpha}_{g}([0, T], \mathbb{R})$ -- a geometric rough path above $B^{H},$ $L = \Delta$ -- negative self-adjoint operator generating an analytic semigroup $\{S_{t}\}_{t\geq0}$ on $\mathcal{U}$ (see \cite[P. 421]{bib15}), 
	$D(\Delta) = \big\{v(\cdot) \in H^{2}([0, 1], \mathbb{R}): v(0) = v(1) = 0\big\},$ $(F(y)z)(x) = y(x)z,$ $y \in\mathcal{U},$ $z \in \mathbb{R},$ $x \in [0, 1].$ Since $F(y)z = yz,$ $y \in\mathcal{U}_{\beta},$ $z \in \mathbb{R},$ the function $F$ and its Frechet derivatives are locally bounded on $\mathcal{U}_{\beta} := D((-\Delta)^{\beta}),$ $\beta \in \mathbb{R}.$
\end{example}

\section*{Declarations}

Authors declare that they have no conflict of interest.

\begin{appendices}

\section{Proofs of auxiliary statements}\label{secA1}

	\subsection{Proof of Proposition \ref{th2.5}} Is is enough to prove the statement for $\mathcal{U} = \ell_{2}.$ We will use a similar idea to \cite[Sec. 5]{bib19}.
	
	Let's consider the Hilbert space $\ell_{2} \times \mathbb{R}$ of sequences
	\begin{center}
	$(x_{0},x_{1},x_{2}, ...),$ $x_{0} \in \mathbb{R}, (x_{1},x_{2}, ...) \in \ell_{2}.$
	\end{center}
	with the inner product
	$$
	(x, y) = \sum_{n = 1}^{\infty}x_{n}y_{n} + x_{0}y_{0} \ \forall x, y \in \ell_{2} \times \mathbb{R}.
	$$
	Let's define the bijective mapping $\varphi$ from  $\ell_{2} \cup \{\mbox{\mathversion{bold}$\infty$}\}$ into the closed subset $\displaystyle \{y \in \ell_{2} \times \mathbb{R}: \sum_{n = 1}^{\infty}y_{n}^{2} + y_{0}^{2} = y_{0}\}$ by the following:
	\begin{center}
	$\displaystyle y_{0} = \frac{\displaystyle \sum_{n = 1}^{\infty}x_{n}^{2}}{\displaystyle 1 + \sum_{n = 1}^{\infty}x_{n}^{2}}, y_{i} = \frac{\displaystyle x_{i}}{\displaystyle 1 + \sum_{n = 1}^{\infty}x_{n}^{2}} \ \forall i \in \mathbb{N},$ if $x = (x_{1},x_{2}, ...) \in \ell_{2}$
	\end{center}
	and
	\begin{center}
	$\displaystyle y_{0} = 1, y_{i} = 0 \ \forall i \in \mathbb{N},$ if $x = \mbox{\mathversion{bold}$\infty$}.$
	\end{center}
	Let's define a metric on $\widehat{\ell}_{2}=\ell_{2} \cup \{\mbox{\mathversion{bold}$\infty$}\}$:
	$$
	\widehat{d}(x, \tilde{x})^{2} = \sum_{n = 1}^{\infty}(y_{n} - \tilde{y}_{n})^{2} + (y_{0} - \tilde{y}_{0})^{2}
	$$
	$$
	= \sum_{n = 1}^{\infty}\Big(\frac{x_{n}}{1 + \|x\|^{2}} - \frac{\tilde{x}_{n}}{ 1 + \|\tilde{x}\|^{2}}\Big)^{2} + \Big(\frac{\|x\|^{2}}{1 + \|x\|^{2}} - \frac{\|\tilde{x}\|^{2}}{1 + \|\tilde{x}\|^{2}}\Big)^{2} 
	$$
	$$
	= \frac{\displaystyle \sum_{n = 1}^{\infty}\big(x_{n}(1 + \|\tilde{x}\|^{2}) - \tilde{x}_{n}(1 + \|x\|^{2})\big)^{2} + \big(\|x\|^{2}(1 + \|\tilde{x}\|^{2}) - \|\tilde{x}\|^{2}(1 + \|x\|^{2})\big)^{2}}{(1 + \|x\|^{2})^{2}(1 + \|\tilde{x}\|^{2})^{2}} 
	$$
	$$
	= \frac{(1 + \|x\|^{2})(1 + \|\tilde{x}\|^{2})\|x - \tilde{x}\|^{2}}{(1 + \|x\|^{2})^{2}(1 + \|\tilde{x}\|^{2})^{2}} = \frac{\|x - \tilde{x}\|^{2}}{(1 + \|x\|^{2})(1 + \|\tilde{x}\|^{2})} \ \forall x, \tilde{x} \in \ell_{2}. 
	$$
	If $x \in \ell_{2},$ $\tilde{x} = \mbox{\mathversion{bold}$\infty$},$ then
	$$
	\widehat{d}(x, \mbox{\mathversion{bold}$\infty$})^{2} = \sum_{n = 1}^{\infty}y_{n}^{2} + (y_{0} - 1)^{2} = (\sum_{n = 1}^{\infty}y_{n}^{2} + y_{0}^{2}) + 1 - 2y_{0} = y_{0} + 1 - 2y_{0} = 1 - y_{0} 
	$$
	$$
	= 1 - \frac{\|x\|^{2}}{1 + \|x\|^{2}} = \frac{1}{1 + \|x\|^{2}}.
	$$

	Let's show that the metric space $\widehat{\ell}_{2}$ is complete with respect to $\widehat{d}.$ Let $(x^{n}) \subset \widehat{\ell}_{2}$ Be a Caushy sequence, then the corresponding sequence $y^{n}=\varphi(x^{n})$ will be a Caushy sequence in $\ell_{2} \times \mathbb{R}$. Since the sphere $\{\tilde{y} \in \ell_{2} \times \mathbb{R}: \sum_{m = 1}^{\infty}\tilde{y}_{m}^{2} + \tilde{y}_{0}^{2} = \tilde{y}_{0}\}$ is a closed subset of $\ell_{2} \times \mathbb{R},$ $(y^{n})$ converges to some element of this sphere. 
	
	Proposition \ref{th2.5} is proved.

\end{appendices}

\bibliography{sn-bibliography}

\end{document}